\documentclass[12pt]{amsart}

\usepackage[english]{babel}
\usepackage{amssymb}
\usepackage{amsmath}

\begin{document}

\baselineskip=17pt

\title{METRIZABLE DH-SPACES WITH A DENSE COMPLETE SUBSET}

\author{SERGEY MEDVEDEV}
\address{South Ural State University, Chelyabinsk, 454080 Russia}
\email{medv@math.susu.ac.ru}

\begin{abstract}
It is proved that for an $h$-homogeneous space $X$ the following conditions are equivalent: 1) $X$ is a densely homogeneous space with a dense complete subspace;  2) $X$ is $\sigma$-discretely controlled.
\end{abstract}



\maketitle

All spaces under discussion are metrizable and strongly zero-dimensional.

In 2005, M.~Hru{\v{s}}{\'{a}}k and B.~Z.~Avil{\'{e}}s \cite[Theorem 2.3]{HrA} proved that every analytic \textsf{CDH}-space is completely Baire. It is therefore of interest to construct a non-complete \textsf{CDH}-space with a dense complete subset. It was realized by R.~Hern\'{a}ndez-Guti\'{e}rrez, M.~Hru\v{s}\'{a}k, and J.~van Mill \cite[Corollary 4.6]{HHrM} and A.~Medini \cite[Theorem 7.3]{AM-14}. They shown that if a homogeneous space $X$ with a dense complete subset can be embedded in the Cantor set $2^\omega$ by a special way, then $X$ is \textsf{CDH}. Moreover, Medini \cite[Theorem 9.4]{AM-14} proved that $X^\omega$ is \textsf{CDH} if $X$ is a countably controlled space. We suggest a criterion for checking the \textsf{CDH}-property for a space without its embedding into the Cantor set. A similar result is obtained for non-separable spaces.  Notice that the following characterization of non-Baire densely homogeneous $h$-homogeneous spaces was given by Medvedev \cite[Theorem 6]{Med15}:

\textbf{Theorem 1.}
\textit{Let $X$ be a non-$\sigma$-discrete $h$-homogeneous space. Then $X$ is a densely homogeneous space of the first category if and only if every $\sigma$-discrete subset of $X$ is a $G_\delta$-set in $X$.}

For all undefined terms and notations see~\cite{Eng}.

We shall say that a space $X$ is \textit{$h$-homogeneous} if $\mathrm{Ind}X =0$ and every non-empty clopen subset of $X$ is homeomorphic to $X$ (this term was first used by Ostrovsky \cite{ost81}). By a \textit{complete} space we mean a completely metrizable space, or an absolute $G_\delta$-space. A separable topological space $X$ is \textit{countable dense homogeneous} (briefly, \textsf{CDH}) if, given any two countable dense subsets $A$ and $B$ of $X$, there is a homeomorphism  $h \colon X  \rightarrow X$ such that $h(A)= B$. Similarly, a space $X$ is \textit{densely homogeneous} (briefly, \textsf{DH}) if, given any two $\sigma$-discrete dense subsets $A$ and $B$  of $X$, there is a homeomorphism  $h \colon X \rightarrow X$ such that $h(A)= B$.

Medini \cite{AM-14} says that a separable space $X$ is \textit{countably controlled} if for every countable $D \subset X$ there exists a Polish subspace $G$ of $X$ such that $D \subset G$. Similarly, we shall say that a space $X$ is \textit{$\sigma$-discretely controlled} if for every $\sigma$-discrete $D \subset X$ there exists a completely metrizable subspace $G$ of $X$ such that $D \subset G$. Clearly, every $\sigma$-discretely controlled space has a dense complete subspace because a metrizable space always contains a dense $\sigma$-discrete subspace.

If $\mathcal{U}$ is a family of sets, then $\bigcup \mathcal{U} = \bigcup \{U \colon U \in \mathcal{U} \}$.

We shall use technique from \cite{Med12} which based on the following lemma.

\textbf{Lemma 1.}\label{le-res}
\textit{Let $X_i$ be a metrizable space with ${\mathrm{Ind}} X_i =0$ and
$F_i$ is a nowhere dense closed subset of $X_i$, where $i \in \{1,
2 \}$. Let $f \colon F_1\rightarrow F_2$ be a homeomorphism.}

\textit{Then there exist a $\sigma$-discrete (in $X_i$) cover $\mathcal{V}_i$ of $X_i \setminus F_i$ by non-empty pairwise disjoint clopen subsets of $X_i$, $i \in \{1, 2 \}$, and a bijection $\psi \colon \mathcal{V}_1 \rightarrow \mathcal{V}_2$ such that for any subsets $D_1 \subseteq \bigcup \mathcal{V}_1$
and $D_2 \subseteq \bigcup \mathcal{V}_2$ and any bijection $g \colon D_1 \rightarrow D_2$ satisfying $g(D_1 \cap V) = D_2 \cap \psi(V)$ for every $V \in \mathcal{V}_1$, the combination $f \triangledown g \colon F_1 \cup D_1 \rightarrow F_2 \cup D_2$ is continuous at each point of $F_1$ and its inverse $(f \triangledown g)^{-1}$ is continuous at each point of $F_2$.}

Under the conditions of Lemma \ref{le-res}, we say \cite{Med_Closed} that:

1) the cover $\mathcal{V}_i$ \textit{forms a residual family} with
respect to $F_i$, where $i \in \{1, 2 \}$,

2) the bijection $\psi \colon \mathcal{V}_1 \rightarrow
\mathcal{V}_2$ is \textit{agreed to} $f$.

\textbf{Remark 1.} A construction similar to Lemma \ref{le-res} was used by B.~Knaster and M.~Reichbach \cite{KnR}, R.~Hern\'{a}ndez-Guti\'{e}rrez, M.~Hru\v{s}\'{a}k, and J.~van Mill \cite{HHrM}, and A.~Medini \cite{AM-14} for separable spaces and by A.~V.~Ostrovsky \cite{ost81} and the author \cite{msu} for non-separable spaces. The last two papers were written in Russian, therefore the proof of Lemma \ref{le-res} was also given in
\cite{Med_Closed}.

Independently Lemma \ref{le-res} was obtained by F. van Engelen
\cite[Theorem 2.3]{EngInf} for weight-homogeneous spaces. In his notation, the triple $\langle \mathcal{V}_1, \mathcal{V}_2, \psi \rangle$ is called a \textit{Knaster-Reichbach cover}, or a \textsf{KR}-cover, for $\langle X_1 \setminus F_1, X_2 \setminus F_2, f \rangle$.

\section*{The main theorem}

\textbf{Theorem 2.}\label{T2}
\textit{Let $X$ be an $h$-homogeneous space. Then $X$ is a densely homogeneous space with a dense complete subspace if and only if $X$ is $\sigma$-discretely controlled.}

\textsc{Proof.} We first verify that if $X$ is a densely homogeneous space, then $X$ is $\sigma$-discretely controlled. Take a $\sigma$-discrete set $D \subset X$ and expand it to a dense $\sigma$-discrete set $D^* \subset X$. Under the conditions of the theorem $X$ has a dense complete subspace $G$; hence there is a dense $\sigma$-discrete set $A \subset G$. Let us choose a homeomorphism $f \colon X \rightarrow X$ such that $f(A) = D^*$. Then $D$ lies entirely within the completely metrizable space $f(G)$.

Now we shall prove that a $\sigma$-discretely controlled space $X$ is densely homogeneous.

To begin, assume additionally that $X$ is a completely metrizable space.

Clearly, every completely metrizable space is $\sigma$-discretely controlled. The statement is trivial if $X$ is a point. By virtue of \cite[Exercise 4.3.H]{Eng}, the theorem holds if $X$ is the Cantor set or the space of irrationals. According to \cite[Theorem 5.4]{msu} (see also \cite[Corollary 2]{Med12}) the Baire space $B(k)$ of weight $k$ is densely homogeneous.  There is no another $h$-homogeneous complete spaces.

Next, consider the main case when $X$ is not completely metrizable.

Fix a metric $\varrho$ on $X$ which is bounded by 1.

Take two $\sigma$-discrete dense subsets $A$ and $B$ of $X$.

According to the Dowker-Hurewicz theorem (see \cite[Theorem 7.3.1]{Eng}) there exists a sequence $\mathcal{U}_{\, 0}, \mathcal{U}_{\, 1}, \ldots$ of discrete
clopen covers of $X$ such that $\mathcal{U}_{\, n+1}$ is a refinement of $\mathcal{U}_{\, n}$ for each $n$ and the family $\mathcal{U} = \{U \in \mathcal{U}_{\, n} \colon n \in \omega \}$ forms a base for $X$. For each $U \in \mathcal{U}$ fix a homeomorphism $\varphi_U \colon X \rightarrow U$. One can verify that the sets $A^* = A \cup \bigcup \{\varphi_U(A) \colon U \in \mathcal{U} \}$ and $B^* = B \cup \bigcup \{\varphi_U(B) \colon U \in \mathcal{U} \}$ are $\sigma$-discrete. Then the set $A^* \cup B^*$ is $\sigma$-discrete. Since $X$ is $\sigma$-discretely controlled, there exists a complete subspace $G \subset X$ containing $A^* \cup B^*$. Clearly,  $G$ is dense in the $h$-homogeneous space $X$. The Stone theorem \cite[Theorem 1]{St1} implies that there is a homeomorphism $\chi \colon B(\tau) \rightarrow G$, where $\tau =w(X)$. Fix the standard metric on the Baire space $B(\tau)$. By construction, $\varphi_U(F) \cap (A \cup B)= \emptyset$ for any $U \in \mathcal{U}$ and $F \subset X \setminus (A^* \cup B^*)$.

Let $A=\bigcup \{A_n \colon n \in \omega\}$ and $B =\bigcup \{B_n \colon n \in \omega\}$, where each $A_n, \, B_n$ is a closed discrete subset of $X$, and $X \setminus G = \bigcup \{X_n \colon n \in \omega \}$, where each $X_n$ is a closed nowhere dense subset of $X$.

We will construct a homeomorphism $f \colon X \rightarrow X$ satisfying $f(A) = B$ by induction on $n$. For each $n \in \omega$ we will define a homeomorphism $f_n \colon X \rightarrow X$, closed nowhere dense sets $F_n$ and $E_n$, a \textsf{KR}-cover $\langle \mathcal{V}_n, \mathcal{W}_n, \psi_n \rangle$ for $\langle X \setminus F_n, X \setminus E_n, f_n|_{F_n} \rangle$ such that:

(a) $f_n( F_n) = E_n$,

(b) $F_n \cup X_n \cup A_n \subset F_{n+1}$,

(c) $E_n \cup X_n  \cup B_n \subset E_{n+1}$,

(d) $f_n(A \cap F_n) = B \cap E_n$,

(e)$f_m(x) = f_n(x)$ for each $x \in F_n$ whenever $m \geq n$,

(f) $\mathcal{V}_{n+1}$ refines $\mathcal{V}_n$ and
$\mathcal{W}_{n+1}$ refines $\mathcal{W}_n$,

(g) if $V \in \mathcal{V}_n$ and $V^* \in \mathcal{V}_m$ for $n < m$,
     then $V^* \subset V$ if and only if $\psi_m(V^*) \subset \psi_n(V)$,

(h) $\mathrm{mesh}(\mathcal{V}_n) \leq 1/(n+1)$ and $\mathrm{mesh}(\mathcal{W}_n) \leq 1/(n+1)$,

(i) $\mathrm{diam}(\chi^{-1}(U \cap G)) \leq 1/(n+1)$ for each $U \in \mathcal{V}_n \cup \mathcal{W}_n$.

For $n=0$, let $F_0 =E_0= \emptyset$, $\mathcal{V}_0 =
\mathcal{W}_0 = \{X \}$, $\psi_0(X)= X$, and $f_0(x)=x$.

Now let us consider the induction step $n+1$ for some $n \in \omega$.

First, fix $V \in \mathcal{V}_n$ and $W = \psi_n(V) \in \mathcal{W}_n$. Let $Y = V \cap X_m$, where $m$ is the least $i$ with $V \cap X_i \neq \emptyset$, and $Z = W \cap X_l$, where $l$ is the least $i$ with $W \cap X_i \neq \emptyset$.

The open set $V \setminus Y$ contains the set $V \cap G$. Then for each point $x \in V \cap G $ we can find its neighborhood $U_x \in \mathcal{U_\ell}$ such that $U_x \subset V \setminus Y$, $\mathrm{diam}(U_x) \leq 1/(n+2)$, $\mathrm{diam}(\chi^{-1}(U_x \cap G)) \leq 1/(n+2)$, and $\ell$ is the least suitable index. The family $\mathcal{V}^*_V= \{U_x \colon x \in V \cap G\}$ consists of pairwise disjoint clopen subsets of $X$, $\mathcal{V}^*_V \subset \mathcal{U}$, and $V \cap G  \subset \bigcup \mathcal{V}^*_V  \subset V \setminus Y$. Clearly, $Y \subset V \setminus \bigcup \mathcal{V}^*_V$. Note that $V \setminus \bigcup \mathcal{V}^*_V$ is a closed nowhere dense subset of $V$ because the set $V \cap G$ is dense in $V$.

Likewise, choose a family $\mathcal{W}^*_V= \{U_x \colon x \in W \cap G \}$ such that $\mathcal{W}^*_V$ consists of pairwise disjoint clopen subsets of $X$, $\mathcal{W}^*_V \subset \mathcal{U}$, $\mathrm{diam}(U) \leq 1/(n+2)$, $\mathrm{diam}(\chi^{-1}(U \cap G)) \leq 1/(n+2)$ for each $U \in \mathcal{W}^*_V$, and $W \cap G \subset \bigcup \mathcal{W}^*_V  \subset  W \setminus Z$.

Since $X$ is locally of weight $\tau$, there exists a closed discrete set $A^\prime \subset A \cap V$ of cardinality $\tau$ such that $A_{m} \cap V \subset A^\prime$, where $m$ is the least $i$ with $V \cap A_i \neq \emptyset$. The union $Y \cup A^\prime$ is a closed nowhere dense subset of $V$. Hence $V$ contains a set $U \in \mathcal{V}^*_V$ such that $U \cap (Y \cup A^\prime) = \emptyset$. Similarly, choose a closed discrete set $B^\prime \subset B \cap W$ of cardinality $\tau$ such that $B_{l} \cap W \subset B^\prime$, where $l$ is the least $i$ with $W \cap B_i \neq \emptyset$. Then $W$ contains a set $O \in \mathcal{W}^*_V$ which misses $Z \cup B^\prime$. Fix a bijection $g \colon A^\prime \rightarrow B^\prime$.

By construction, $A^\prime$, $V \setminus \bigcup \mathcal{V}^*_V$, and $ \varphi_U(W \setminus \bigcup \mathcal{W}^*_V )$ are disjoint closed nowhere dense subsets of $V$. Similarly, $B^\prime$, $W \setminus \bigcup \mathcal{W}^*_V$, and $\varphi_{O}(V \setminus \bigcup \mathcal{V}^*_V)$ are disjoint closed nowhere dense subsets of $W$. Put
\[\begin{array}{c}
F_V = A^\prime \cup (V \setminus \bigcup \mathcal{V}^*_V) \cup \varphi_U(W \setminus \bigcup \mathcal{W}^*_V ),  \\
E_V = B^\prime \cup (W \setminus \bigcup \mathcal{W}^*_V ) \cup \varphi_{O}(V \setminus \bigcup \mathcal{V}^*_V).
\end{array} \]

Notice that $A \cap F_V = A^\prime$ and $B \cap E_V = B^\prime$.

Define the homeomorphism $f_V \colon F_V \rightarrow E_V$ by setting
\[ f_V(x) =\left \{
\begin{array}{cl}
g(x) & \mbox{if  } x \in A^\prime,  \\
\varphi_O(x) & \mbox{if  } x \in V \setminus \bigcup \mathcal{V}^*_V, \\
(\varphi_U)^{-1}(x) & \mbox{if  } x \in \varphi_U(W \setminus \bigcup \mathcal{W}^*_V ).
\end{array}   \right. \]

By Lemma \ref{le-res}, there exists a \textsf{KR}-cover $\langle \mathcal{V}_V,
\mathcal{W}_V, \psi_V \rangle$ for $\langle V \setminus F_V, W \setminus E_V, f_V \rangle$.  By construction, $\bigcup \mathcal{V}_V = \bigcup \mathcal{V}^*_V \setminus A^\prime$ and $\bigcup \mathcal{W}_V = \bigcup \mathcal{W}^*_V \setminus B^\prime$. By reducing $\mathcal{V}_V $ and $\mathcal{W}_V $ to fragments if it is necessary, we can obtain that $\mathrm{diam}(\chi^{-1}(U \cap G)) \leq 1/(n+2)$ for each $U \in \mathcal{V}_V \cup \mathcal{W}_V$, $\mathcal{V}_V \subset \mathcal{U}$, and $\mathcal{W}_V \subset \mathcal{U}$.

Repeat this construction for each $V \in \mathcal{V}_n$. It is not hard to check that $F_{n+1}  = F_n \cup \bigcup \{F_V \colon V \in \mathcal{V}_n \}$ and $E_{n+1} = E_n \cup \bigcup \{E_V \colon V \in \mathcal{V}_n \}$ are closed nowhere dense subsets of $X$. From Lemma \ref{le-res} it follows that the mapping \[ f_{n+1}(x) =\left \{
\begin{array}{cl}
f_V(x) & \mbox{if  } x \in F_V \mbox{ for some } V
\in \mathcal{V}_n,  \\
f_n(x) & \mbox{if  } x \in F_n
\end{array}   \right. \]
is a homeomorphism between $F_{n+1}$ and $E_{n+1}$.

Define $\mathcal{V}_{n+1} = \{U \in \mathcal{V}_V \colon V \in \mathcal{V}_{n} \}$ and $\mathcal{W}_{n+1} = \{\mathcal{W}_V \colon V \in \mathcal{V}_{n} \}$. Then $\psi_{n+1}= \bigcup \{\psi_V \colon V \in \mathcal{V}_n \}$ is a bijection between $\mathcal{V}_{n+1}$ and $\mathcal{W}_{n+1}$ which is agreed to $f_{n+1}$. Thus, $\langle \mathcal{V}_{n+1}, \mathcal{W}_{n+1}, \psi_{n+1} \rangle$ is a \textsf{KR}-cover for $\langle X \setminus F_{n+1}, X \setminus E_{n+1}, f_{n+1} \rangle$.

It remains to extend the mapping $f_{n+1}$ from $F_{n+1}$ to the space $X$. Both sets $U \in \mathcal{V}_{n+1}$ and $\psi_{n+1}(U) \in \mathcal{W}_{n+1}$ are clopen in the $h$-homogeneous space $X$. Choose a homeomorphism $g_U \colon U \rightarrow \psi_{n+1}(U)$ and put $f_{n+1}(x)= g_U(x)$ if $x \in U \in \mathcal{V}_{n+1}$. From Lemma \ref{le-res} it follows that $f_{n+1} \colon X \rightarrow X$ is a homeomorphism.

One can check that all conditions (a)--(i) are satisfied.
This completes the induction step.

Let $F= \bigcup \{F_n \colon n \in \omega \}$ and $E= \bigcup \{E_n \colon n \in \omega \}$. If $x \in F_n$ for some $n \in \omega$, set $f(x) =f_n(x)$.  Condition (e) implies that this rule define a bijection between $F$ and $E$. Moreover, from conditions (c) and (d) it follows that $f(A)= B$. For every point $x \in X \setminus F \subset G$ and every $n \in \omega$ there exists the unique clopen set $U^x_n \in \mathcal{V}_n$ that contains $x$.  Clearly, $x = \bigcap \{U^x_n \colon n \in \omega \}$. By virtue of the Cantor theorem \cite[Theorem 4.3.8]{Eng} and condition (i), the intersection $\bigcap \{\chi^{-1} \left(G \cap \psi_n(U^x_n) \right) \colon n \in \omega \}$ is non-empty and consists of a point $z \in B(\tau)$. Condition (h) implies that the intersection $ \bigcap \{\psi_n(U^x_n) \colon n \in \omega \} $ coincides with the point $\chi(z)$. Clearly, $\chi(z) \in X \setminus E$. Set $f(x)= \chi(z)$. Likewise, every point $y \in X \setminus E$ can be represented as $\bigcap \{U^y_n \colon n \in \omega \}$, where each $U^y_n \in \mathcal{W}_n$. Then $y= f(x_0)$ for $x_0 = \bigcap \{\psi_n^{-1}(U^y_n) \colon n \in \omega \} \in X \setminus F$. Hence, we obtain the mapping $f \colon X \rightarrow X$ such that $f(X)=X$. One easily verifies that $f$ is a bijection.

From the definition of a \textsf{KR}-cover it follows that $f$ is continuous on $F$ and $f^{-1}$ is continuous on $E$. Condition (h) implies that the sequence $\{U^x_n \colon n \in \omega \}$ is a base at the point $x \in X \setminus F$ and the sequence $\{U^y_n \colon n \in \omega \}$ is a base at the point $y \in X \setminus E$. By construction, if $y = f(x)$ then $f(U^x_n)=\varphi_n(U^x_n) = U^y_n$. From this it follows that $f$ is continuous on $X \setminus F$ and $f^{-1}$ is continuous on $X \setminus E$. Hence, $f$ is a homeomorphism.
$\square$

\textbf{Corollary 1.}
\textit{Every countably controlled  $h$-homogeneous space is \textsf{CDH}.}

Let us show that an $h$-homogeneous $\sigma$-discretely controlled space satisfies a strengthening version of dense homogeneity.

\textbf{Theorem 3.}
\textit{Let $X$ be an $h$-homogeneous $\sigma$-discretely controlled space. Suppose that $\{A_n \colon n \in \omega \}$ and $\{B_n \colon n \in \omega \}$ are two families of $\sigma$-discrete dense subsets of $X$ such that $A_i \cap A_j = \emptyset$ and $B_i \cap B_j = \emptyset$ whenever $i \neq j$. Then there is a homeomorphism $f \colon X \rightarrow X$ such that $f(A_n) = B_n$ for each $n \in \omega$.}

\textsc{Proof. }The proof is similar to the corresponding part of the proof of Theorem \ref{T2}. We only point out the changes and omit like details. Use the notation from the proof of Theorem \ref{T2}.

For each $n \in \omega$ let $A_n=\bigcup \{A^i_n \colon i \in \omega\}$ and $B_n =\bigcup \{B^i_n \colon i \in \omega\}$, where $A^i_n$ and $ B^i_n$ are closed discrete subsets of $X$. Define $A = \bigcup \{A_n \colon n \in \omega \}$ and $B = \bigcup \{B_n \colon n \in \omega \}$. Put
\[A^* \cup B^* = (A \cup B) \cup \bigcup \{\varphi_U(A \cup B) \colon U \in \mathcal{U} \}.\]
Fix a complete subspace $G \subset X$ containing the $\sigma$-discrete set $A^* \cup B^*$.

For each $n \in \omega$ we will construct a homeomorphism $f_n \colon X \rightarrow X$ satisfying conditions (a)--(i), but conditions (b)--(d) will be replaced by the following ones:

(b*) $F_n \cup X_n \cup \bigcup \{A^i_n \colon i \leq n \} \subset F_{n+1}$,

(c*) $E_n \cup X_n  \cup \bigcup \{B^i_n \colon i \leq n \} \subset E_{n+1}$.

(d*) $f_n(A_j \cap F_n) = B_j \cap E_n$ for each $j \leq n$.

For this, given $V \in \mathcal{V}_n$, $W = \psi_n(V) \in \mathcal{W}_n$, and $j \leq n$, choose closed discrete sets $A^\prime(j) \subset A_j \cap V$ and $B^\prime(j) \subset B_j \cap W$ each of cardinality $\tau$ such that $(A^i_j \cap V) \setminus F_n \subset A^\prime(j)$ and $(B^i_j \cap W) \setminus E_n \subset B^\prime(j)$ for every $i \leq n$. Fix a bijection $g_j \colon A^\prime(j) \rightarrow B^\prime(j)$. Under the conditions of the theorem, $A^\prime(j_1) \cap A^\prime(j_2) = \emptyset$ and $B^\prime(j_1) \cap B^\prime(j_2) = \emptyset$ whenever $j_1 \neq j_2$. We thus get a bijection $g \colon A^\prime \rightarrow B^\prime$, where $A^\prime  =  \bigcup \{A^\prime(j) \colon j \leq n \}$ and  $B^\prime  =  \bigcup \{B^\prime(j) \colon j \leq n \}$. We omit further details of the proof. $\square$

\end{document}